\documentclass[12pt]{article}
\usepackage{amssymb,latexsym,amsmath,euscript}\usepackage{amsfonts}
\usepackage[cp1251]{inputenc}\usepackage[epsf]{}\usepackage[english,russian]{babel}\usepackage{graphicx}
\newtheorem{prop}{Proposition}\newtheorem{cor}{Corollary}
\newtheorem{lm}{Lemma}\newtheorem{thm}{Theorem}
\newtheorem{rema}{Remark}
\newtheorem{defi}{Definition}

\begin{document}

\title{Morse-Smale systems without heteroclinic submanifolds on codimension one separatrices
\footnote{2000{\it Mathematics Subject Classification}. Primary 37D15; Secondary 58C30}
\footnote{{\it Key words and phrases}: Morse-Smale systems, heteroclinic intersections}}

\author{V.~Grines$^{1}$\and V.~Medvedev$^{1}$\and E.~Zhuzhoma$^{1}$}
\date{}
\maketitle

\begin{center}
{\small $^{1}$ National Research University Higher School of Economics,\\ 25/12 Bolshaya Pecherskaya, 603005, Nizhni Novgorod, Russia}
\end{center}

\renewcommand{\abstractname}{Absrtact}\renewcommand{\refname}{Bibliography}
\renewcommand{\figurename}{Fig.}

\begin{abstract}
We study a topological structure of  a closed $n$-manifold $M^n$ ($n\geq 3$) which admits a Morse-Smale  diffeomorphism 
such that codimension one separatrices of saddles periodic points   have no heteroclinic intersections  different from  heteroclinic points. Also we consider gradient like flow on $M^n$
such that   codimension one separatices  of saddle singularities have no intersection at all.   We show that $M^n$ is either an $n$-sphere $S^n$, or the connected sum of a finite number of copies of $S^{n-1}\otimes S^1$ and a finite number of special manifolds $N^n_i$ admitting polar Morse-Smale systems. Moreover, if some $N^n_i$ contains a single saddle, then $N^n_i$ is projective-like (in particular, $n\in\{4,8,16\}$, and $N^n_i$ is a simply-connected and orientable manifold). Given input dynamical data, one constructs a supporting manifold $M^n$.
We give a formula relating the number of sinks, sources and saddle periodic points to the connected sum for $M^n$. As a consequence, we obtain conditions for the existence of heteroclinic intersections for Morse-Smale diffeomorphisms and a periodic trajectory for Morse-Smale flows.
\end{abstract}

\section*{Introduction}\label{s:title-introd}\nopagebreak

Morse-Smale dynamical systems (flows and diffeomorphisms) was introduced by Steve Smale \cite{Smale60a} in 1960 using the explicit description of properties of a non-wandering set. At his work, S.~Smale grasped an intimate connection between dynamic properties of Morse-Smale systems and topological properties of ambient manifolds. He discovered in a spirit of Morse's inequalities \cite{Morse-book-1934} the relations between numbers of periodic points and the Betty numbers $\beta_0(M^n)$, $\beta_1(M^n)$, $\ldots$, $\beta_n(M^n)$ of the ambient $n$-manifolds $M^n$, where $\beta_i(M^n)=rank~H_i(M^n,\mathbb{Z})$. After paper \cite{Smale60a}, many results was obtained for various classes of Morse-Smale systems concerning relations between dynamics and the topological structure of ambient manifolds, see the book \cite{grinmedpoch-book-2016-eng} with numerous references. In the present paper  we study the topological structure of closed manifolds supporting 
a Morse-Smale  diffeomorphism or a gradient like flow such that codimension one separatrices of saddles periodic points of diffeomorphisms  have no heteroclinic intersections  different from  heteroclinic points and  codimension one separatrices of saddle singularities of flows have no intersection at all. The  paper is generalization and specification of results which was obtained by the authors in \cite{GrinesMedvedevZh2017}. 

We identify a continuous-time dynamical system with a flow, while a discrete-time dynamical system is identified with a diffeomorphism.
Many definitions for diffeomorphisms and flows are similar. So, we shall give mainly the notation for diffeomorphisms giving the exact notation for flows if necessary. Let $f: M^n\to M^n$ be a diffeomorphism of $n$-manifold $M^n$, and $p$ a periodic point of period $k\in\mathbb{N}$. The stable manifold $W^s(p)$ is defined to be the set of points $x\in M^n$ such that $\varrho(f^{kj}(x);p)\to 0$ as $j\to\infty$ where $\varrho$ is a metric on $M^n$. The unstable manifold $W^u(p)$ is the stable manifold of $p$ for the diffeomorphism $f^{-1}$. Stable and unstable manifolds are called invariant manifolds. It is well known that if $p$ is hyperbolic, then every invariant manifold is an immersed submanifold homeomorphic to Euclidean space. Moreover, $W^s(p)$ and $W^u(p)$ are intersected transversally at $p$, and $\dim W^s(p)+\dim W^u(p)=n$.

A dynamical system $\mathcal{F}$  is \textit{Morse-Smale} if it is structurally stable and the non-wandering set $NW(\mathcal{F})$ of $\mathcal{F}$ consists of a finitely many periodic orbits\footnote{We consider fixed points being periodic orbits of the trivial period 1.} if $\mathcal{F}$ is diffeomorphism  and a finitely many equilibrium states and closed trajectories if $\mathcal{F}$ is flow.   From the modern point of view, the Morse-Smale systems are exactly the structurally stable dynamical systems with zero topological entropy. Morse-Smale flow is called \textit{gradient-like} if its nonwandering set does not contain closed trajectories.    A periodic orbit $p$  is called a \textit{sink} (resp. \textit{source}) one if $\dim W^s(p)=n$ and $\dim W^u(p)=0$ (resp. $\dim W^s(p)=0$ and $\dim W^u(p)=n$). A sink or source periodic orbit is called a \textit{node} periodic orbit.
A periodic point $\sigma$ is called a \textit{saddle} one if $1\leq\dim W^u(\sigma)\leq n-1$,  $1\leq\dim W^s(\sigma)\leq n-1$.
A component of $W^u(\sigma)\setminus\sigma$ denoted by $Sep^u(\sigma)$ is called an
\textit{unstable separatrix} of $\sigma$. If $\dim W^u(\sigma)\geq 2$, then $Sep^u(\sigma)$ is unique. The similar notation holds for a stable separatrix.
A saddle periodic point $\sigma$ is called \textit{codimension one} if either $\dim W^u(\sigma)=1$, $\dim W^s(\sigma)=n-1$ or
$\dim W^u(\sigma)=n-1$, $\dim W^s(\sigma)=1$. Another words, one of the separatrices, $Sep^s(\sigma)$ in the first case and $Sep^u(\sigma)$ in the second one, is codimension one.  Similar is defined codimension one  saddle equilibrium state  for Morse-Smale flow. 

For any different saddle periodic points (equilibrium states) $p$, $q$  of Morse-Smale diffeomorphism (of gradient like flow)   intersection $W^u(p)\cap W^s(q)$ is either  empty or  transversal.
If  $W^u(p)\cap W^s(q)\neq\emptyset$  the intersection $W^u(p)\cap W^s(q)$ is called \textit{heteroclinic}. Due to the transversality $W^u(p)\pitchfork W^s(q)$, a heteroclinic intersection is either union of coutable number of  isolated points (called \textit{heteroclinic points})  or union of disjoint  $m$-dimensional submanifolds with $m\geq 1$ (called \textit{heteroclinic submanifolds}). 

\begin{defi}We say that a Morse-Smale diffeomorphism  (gradient-like flow) $f: M^n\to M^n$  ($f^t$ on $M^n$) is \textit{without heteroclinic submanifolds on codimension one separatrices} if given any codimension one saddle periodic point (saddle equilibrium state) $p\in NW(f)$,  ($p\in NW(f^t)$)  the codimension one separatrix of $p$ does not contain heteroclinic submanifolds\footnote{If $f: M^n\to M^n$ is diffeomorphism we allow  that  codimension one separatrices can admit heteroclinic points.}.
\end{defi}
Recall that a Morse-Smale diffeomorphism or gradient-like flow  $\mathcal{F}$ is called  \textit{polar} if the non-wandering set $NW(\mathcal{F})$ contains exactly one sink periodic orbit and one source periodic orbit.

Let $\mathbb{S}^n$ be the canonical $n$-sphere defined by the quality $x_1^2+\cdots+x_{n+1}^2=1$ in $\mathbb{R}^{n+1}$ endowed with the coordinates $(x_1;\ldots;x_{n+1})$. By $S^n$ we mean an $n$-sphere homeomorphic to $\mathbb{S}^n$. The open ball $\mathbb{B}^{n+1}$ is defined by the inequality $x_1^2+\cdots+x_{n+1}^2<1$. Denote by $S^{n-1}\otimes S^1$ a total manifold of locally trivial fiber bundle with the base a circle $S^1$ and the fiber a sphere $S^{n-1}$. The manifold $S^{n-1}\otimes S^1$ can be obtained as follows. Take the prime product $S^{n-1}\times [0;1]$ with the natural projection $\pi : S^{n-1}\times [0;1]\to [0;1]$. Gluing the spheres $S^{n-1}\times\{0\}$, $S^{n-1}\times\{1\}$ by a homeomorphism
$S^{n-1}\times\{0\} \to S^{n-1}\times\{1\}$, one gets the total space $S^{n-1}\otimes S^1$ with the projection $S^{n-1}\otimes S^1\to S^1$ induced by $\pi$.  Since any two homeomorphisms $S^{n-1}\to S^{n-1}$ are isotopic provided they either preserve orientation or reverse orientation $S^{n-1}$ sumutaniously, there are only two type of the bundles $S^{n-1}\otimes S^1$: globally non-trivial (skew), and trivial (prime product) $S^{n-1}\times S^1$.

The following theorem describes the topological structure of a supporting manifold admitting Morse-Smale diffeomorphism without hetroclinic submanifolds on codimesion one separatrices. Let us emphasyse that in the paper \cite{GrinesMedvedevZh2017} the  similar result was obtained in suggestions of absence of heteroclinic intersection on   separatrices of any dimension.

\begin{thm}\label{thm:rough-decomposition}
Let $M^n$ be a closed $n$-manifold, $n\ge 3$, supporting a Morse-Smale diffeomorphism $f$ without heteroclinic submanifolds on codimension one separatrices. Suppose that the non-wandering set $NW(f)$ consists of $\mu$ node periodic points, $\nu$ codimension one saddle periodic points, and arbitrary number of saddle periodic points that are not codimension one.
Then the number
 $$ g=\frac{1}{2}\left(\nu - \mu +2\right)\geq 0 $$
is integer. In addition,

1) if $g=0$ then $M^n$ is either $S^n$ or
\begin{equation}\label{eq:decomposition-1}
    M^n=N^n_1\sharp\cdots\sharp\, N^n_{l}
\end{equation}
for some $1\leq l\leq 1+\nu$,  where every $N^n_i$ admits a polar Morse-Smale diffeomorphism without codimension one saddle periodic orbit. Moreover, if $M^n$ is orientable
then $\beta_1(N^n_i)=rank~H_1(M^n,\mathbb{Z})=\beta_{n-1}(N^n_i)=rank~H_{n-1}(M^n,\mathbb{Z})=0$ for $i=1,\ldots,l$;

2) if $g>0$ then $M^n$ is either
\begin{equation}\label{eq:decomposition-2}
    M^n=\underbrace{\left(S^{n-1}\otimes S^1\right)\sharp\cdots\sharp\left(S^{n-1}\otimes S^1\right)}_{g}
\end{equation}
or
\begin{equation}\label{eq:decomposition-3}
    M^n=\underbrace{\left(S^{n-1}\otimes S^1\right)\sharp\cdots\sharp\left(S^{n-1}\otimes S^1\right)}_{g}
    \sharp\, N^n_1\sharp\cdots\sharp\, N^n_{l}
\end{equation}
for some $1\leq l\leq k=\frac{1}{2}\left(\mu + \nu\right)$ where every $N^n_i$ admits a polar Morse-Smale diffeomorphism without codimension one saddle periodic orbits.
Moreover, if $M^n$ is orientable then  $\beta_1(N^n_i)=\beta_{n-1}(N^n_i)=0$ for $i=1,\ldots,l$.
\end{thm}


\begin{cor}\label{cor:for-flows-from-thm:rough-decomposition}
Let $M^n$ be a closed $n$-manifold, $n\ge 3$, supporting a gradient-like  flow $f^t$ without heteroclinic manifolds on codimension one separatrices. Suppose that the non-wandering set $NW(f^t)$ of $f^t$ consists of $\mu$ node equilibrium states, $\nu$ codimension one saddle equilibrium states, and arbitrary number of saddle equilibrium states that are not codimension one.
Then the number
 $$ g=\frac{1}{2}\left(\nu - \mu +2\right)\geq 0 $$
is integer. In addition,

1) if $g=0$ then $M^n$ is either $S^n$ or
\begin{equation}\label{eq:flow-decomposition-1}
    M^n=N^n_1\sharp\cdots\sharp N^n_{l}
\end{equation}
for some $1\leq l\leq 1+\nu$,  where every $N^n_i$ admits a polar Morse-Smale flow without periodic trajectories and codimension one saddle equilibrium states
Moreover, $\beta_1(N^n_i)=\beta_{n-1}(N^n_i)=0$ for $i=1,\ldots,l$ provided $M^n$ is orientable.

2) if $g>0$ then $M^n$ is either
\begin{equation}\label{eq:flow-decomposition-2}
    M^n=\underbrace{\left(S^{n-1}\otimes S^1\right)\sharp\cdots\sharp\left(S^{n-1}\otimes S^1\right)}_{g}
\end{equation}
or
\begin{equation}\label{eq:flow-decomposition-3}
    M^n=\underbrace{\left(S^{n-1}\otimes S^1\right)\sharp\cdots\sharp\left(S^{n-1}\otimes S^1\right)}_{g}
    \sharp N^n_1\sharp\cdots\sharp N^n_{l}
\end{equation}
for some $1\leq l\leq k=\frac{1}{2}\left(\mu + \nu\right)$ where every $N^n_i$ admits a polar Morse-Smale flow without periodic trajectories and codimension one saddle equilibrium states
Moreover, if  $M^n$ is orientable then $\beta_1(N^n_i)=\beta_{n-1}(N^n_i)=0$ for $i=1,\ldots,l$.
\end{cor}
\textsl{Proof of Corollary \ref{cor:for-flows-from-thm:rough-decomposition}}.
If a Morse-Smale flow  $f^t$ has no periodic trajectories, then the time one shift along the trajectories $f^t$ is a Morse-Smale diffeomorphism. Thus the result follows from Theorem \ref{thm:rough-decomposition}.
$\Box$

\begin{rema} For $n=3$, Theorem \ref{thm:rough-decomposition} is a generalization of the main result of \cite{BonattiGrinesMedvedevPecou2002b}.
\end{rema}

To formulate the next result let us introduce the definition of projective-like manifold that is emphasized by the following description of a projective plane. Let $\mathbb{B}^2$ be an open 2-ball with the boundary $S^1=\partial\mathbb{B}^2$ a circle. The identification of opposite points $(a_1;a_2)$, $(-a_1;-a_2)$ of $S^1$ gives the factor-space $S^1/(a_1;a_2)\sim (-a_1;-a_2)$ that is homeomorphic to $S^1$. The natural projection
 $$ S^1\to S^1/\left((a_1;a_2)\sim (-a_1;-a_2)\right)\quad = \,\, S^1 $$
gives the locally trivial fiber bundle $(\partial\mathbb{B}^2=S^1,S^1,S^0)$ with the fiber a zero-dimensional circle $S^0$ that is the union of two points. The projective plane $\mathbb{P}^2$ is obtained from the closed 2-ball $clos~\mathbb{B}^2=\mathbb{B}^2\cup\partial\mathbb{B}^2$ by the identification of every fiber of the locally trivial fiber bundle $(\partial\mathbb{B}^2=S^1,S^1,S^0)$ with a point.

Consider the well-known Hoph fiber bundle $(S^{2n-1},S^n, S^{n-1})$ where $n\in\{2,4,8\}$ \cite{EellsKuiper62,Novikov-book-1976}.
The total space is the $(2n-1)$-sphere $S^{2n-1}$ that projects to the base $S^n$ under the projection $p: S^{2n-1}\to S^n$, and the fiber $p^{-1}(m)$ homeomorphic to $S^{n-1}$ for any
$m\in S^n$.
Take the open balls $\mathbb{B}^4$, $\mathbb{B}^6$, $\mathbb{B}^{16}$ of the dimensions 4, 8, and 16 respectively. Let $N^{2n}$, $n\in\{2,4,8\}$, be a set obtained from the closed $2n$-ball $clos~\mathbb{B}^{2n}=\mathbb{B}^{2n}\cup\partial\mathbb{B}^{2n}$ after the identification of every fiber of the Hoph fiber bundle $(\partial\mathbb{B}^{2n}=S^{2n-1},S^n, S^{n-1})$ with a point. Since the fiber bundle $(\partial\mathbb{B}^{2n}=S^{2n-1},S^n, S^{n-1})$ is locally trivial, $N^{2n}$ is endowed by the structure of a closed (topological) manifold. Such $N^{2n}$, $n\in\{2,4,8\}$ (and every manifold homeomorphic to $N^{2n}$) is called a \textit{projective-like manifold}.

\begin{thm}\label{thm:new-proj-decomposition-for-diff}
Let $M^n$ be a closed $n$-manifold, $n\ge 3$, supporting a Morse-Smale diffeomorphism $f$ without heteroclinic submanifolds on codimension one separatrices. Suppose that the non-wandering set $NW(f)$ consists of $\mu$ node periodic points, $\nu$ codimension one saddle periodic points, and a single saddle fixed point that is not codimension one. Then

1) $n\in\{4,8,16\}$;

2) $M^n$ is either
 $$ M^n=\underbrace{\left(S^{n-1}\otimes S^1\right)\sharp\cdots\sharp\left(S^{n-1}\otimes S^1\right)}_{g}\sharp N^n\,\, \mbox{ provided }\,\, g=\frac{1}{2}\left(\nu - \mu +2\right)>0, $$
or
 $$ M^n=N^n\,\, \mbox{ provided }\,\, \nu=\mu-2 $$
where
\begin{itemize}
  \item $N^n$ is the disjoint union of an open ball $B^n$ and an $\frac{n}{2}$-sphere $S^{\frac{n}{2}}$, $N^n=B^n\cup S^{\frac{n}{2}}$, $B^n\cap S^{\frac{n}{2}}=\emptyset$
  \item for $n\in\{8,16\}$, $N^n$ is a projective-like manifold
  \item the homotopy groups $\pi_1(N^n)=\cdots=\pi_{\frac{n}{2}-1}(N^n)=0$, and hence, $N^n$ is simply connected and orientable.
\end{itemize}
\end{thm}

\begin{cor}\label{cor:from-thm:proj-decomposition-and-diff}
Assume that the conditions of Theorem \ref{thm:rough-decomposition} hold. Suppose that there is $N^n_{i_*}$ in the decompositions (\ref{eq:decomposition-1}) or (\ref{eq:decomposition-3}) such that $N^n_{i_*}$ contains exactly one saddle fixed point. Then
\begin{itemize}
  \item $n\in\{4,8,16\}$
  \item $N^n_{i_*}$ is the disjoint union of an open ball $B^n$ and an $\frac{n}{2}$-sphere $S^{\frac{n}{2}}$, $N^n_{i_*}=B^n\cup S^{\frac{n}{2}}$, $B^n\cap S^{\frac{n}{2}}=\emptyset$
  \item for $n\in\{8,16\}$, $N^n_{i_*}$ is a projective-like manifold
  \item the homotopy groups $\pi_1(N^n_{i_*})=\cdots=\pi_{\frac{n}{2}-1}(N^n_{i_*})=0$, and hence, $N^n_{i_*}$ is simply connected and orientable.
\end{itemize}
\end{cor}

\begin{thm}\label{thm:new-proj-decomposition-for-flow}
Let $M^n$ be a closed $n$-manifold, $n\ge 3$, supporting a Morse-Smale flow $f^t$ without heteroclinic submanifolds on codimension one separatrices. Suppose that the non-wandering set $NW(f^t)$ consists of $\mu$ node equilibrium states, $\nu$ codimension one saddle equilibrium states, and a single saddle equilibrium state that is not codimension one.
Then

1) $n\in\{4,8,16\}$;

2) $M^n$ is either
 $$ M^n=\underbrace{\left(S^{n-1}\otimes S^1\right)\sharp\cdots\sharp\left(S^{n-1}\otimes S^1\right)}_{g}\sharp N^n\,\, \mbox{ provided }\,\, g=\frac{1}{2}\left(\nu - \mu +2\right)>0, $$
or
 $$ M^n=N^n\,\, \mbox{ provided }\,\, \nu=\mu-2 $$
where
\begin{itemize}
  \item $N^n$ is a projective-like manifold
  \item the homotopy groups $\pi_1(N^n)=\cdots=\pi_{\frac{n}{2}-1}(N^n)=0$, and hence, $N^n$ is simply connected and orientable.
\end{itemize}
\end{thm}

The existence of manifold $N^n$ containing exactly one saddle fixed point that is not codimension one follows from \cite{EellsKuiper62}, where one proved the existence
of closed manifolds admitting Morse functions with exactly three critical points, and \cite{Smale60a} where one proved that any gradient flow can be approximated
by Morse-Smale gradient flow.

After Theorem \ref{thm:rough-decomposition}, it is natural to consider the question of so-called realization. Note that the conditions of Theorem \ref{thm:rough-decomposition} have no data on saddles that is not codimension one. Therefore, the summands $N^n_1\sharp\cdots\sharp N^n_{l}$ absence in the following statement. Similar result holds for Morse-Smale flows.
\begin{thm}\label{thm:realization}
Given any integers $\nu\geq 0$ and $\mu\geq 2$ satisfying the following conditions
 $$ \nu\geq\mu-2,\quad g=\frac{1}{2}\left(\nu - \mu +2\right)\in\mathbb{Z}, $$
there is a Morse-Smale diffeomorphism $f: M^n\to M^n$ without heteroclinic submanifolds on codimension one separatrices such that the non-wandering set $NW(f)$ consists of $\mu$ node periodic points, $\nu$ codimension one saddle periodic points where $M^n$ is either
$$ M^n=S^n \quad \mbox{ provided }\quad g=0, $$
or
$$ M^n=\underbrace{\left(S^{n-1}\otimes S^1\right)\sharp\cdots\sharp\left(S^{n-1}\otimes S^1\right)}_{g}\quad \mbox{ provided }\quad g\geq 1. $$
\end{thm}

Let us represent two applications.
The application of Theorem \ref{thm:rough-decomposition} concerns the existence of heteroclinic intersections of codi\-men\-sion one separatrices that form codimension two submanifolds. In the particular case $n=3$, the heteroclinic intersections of two-dimensional separatrices consist of heteroclinic curves. Note that heteroclinic curves is often the mathematical model of so-called separators considered in Solar Magnetohydrodynamics \cite{GrGurZhKaZinina2014,GrinesMedvedevPochZh2015a,Priest-book-1985}. From the modern point of view, reconnections of Solar magnetic field along separators are responsible for Solar flairs \cite{Longcope1996,MacleanBeveridgeLongcopeBrownPriest2005,Priest-Forbs-book-2006}.

\begin{cor}\label{cor:cor1-conditions-for-existing-heter-curves}
Let $f: M^n\to M^n$ be a preserving orientation Morse-Smale diffeomorphism of closed orientable $n$-manifold $M^n$, $n\geq 3$. Suppose that the non-wandering set $NW(f)$ consists of $\mu$ node periodic points, $\nu$ codimension one saddle periodic points, and arbitrary number of saddle periodic points that are not codimension one. If $g=\frac{\nu  - \mu + 2}{2}\geq 1$ and the fundamental group $\pi_1(M^n)$ doesn't contain a subgroup isomorphic to the free product
 $$ \underbrace{\mathbb{Z}*\cdots *\mathbb{Z}}_{g} $$
then there exists saddle periodic points $p$, $q\in NW(f)$ such that $W^s(p)\cap W^u(q)\neq\emptyset$ and the Morse index of $p$ equals $1$, and the Morse index of $q$ equals $n-1$.
\end{cor}

The application of Corollary \ref{cor:for-flows-from-thm:rough-decomposition} concern the existence of a periodic trajectory for Morse-Smale flows.
\begin{cor}\label{cor:periodic-trajectories -exist-if}
Let $f^t$ be a Morse-Smale flow without heteroclinic intersections on a closed orientable manifold $M^n$ of dimension $n\ge 3$, and assume that the non-wandering set $NW(f^t)$ contains the following list of states of equilibrium : $\mu$ nodes, $\nu$ codimension one saddles, and arbitrary number of saddles that are not codimension one.
If $g = \frac{1}{2}\left(\nu - \mu +2\right)\geq 1$ and the fundamental group $\pi_1(M^n)$ doesn't contain a subgroup isomorphic to the free product
 $$ \underbrace{\mathbb{Z}*\cdots *\mathbb{Z}}_{g} $$
then $f^t$ has a periodic trajectory.
\end{cor}

The structure of the paper is the following. In Section \ref{s:def-eng}, we formulate the main definitions and give some previous results. In Section \ref{s:proofs},
we prove main results and its applications.

\textsl{Acknowledgments}.
This work is supported by the Russian Science Foundation under grant 17-11-01041, except the proof of  theorem  \ref{thm:realization}. The latter proof was obtained within the HSE Basic Research Program (project no. 95) in 2018.

\section{Definitions and previous results}\label{s:def-eng}

Here, we recall basic definitions and formulate some results which we need later on. For simplicity, we are considering a discrete-time dynamical system being a diffeomorphism.
For the reference, we formulate the following statement proved in \cite{GrinesMedvedevZh2003a} (see also \cite{GrinesMedvedevPochinkaZh2010,GrinesMedvedevZh2003c}).
\begin{prop}\label{prop:GrinZhMedv2003b-top-closure-of-separatrices}
Let $f: M^n\to M^n$ be a Morse-Smale diffeomorphism, and $Sep^{\tau}(\sigma)$ a separatrix of dimension $1\leq d\leq n-1$ of a saddle fixed point $\sigma$. Suppose that $Sep^{\tau}(\sigma)$ has no intersections with other separatrices. Then $Sep^{\tau}(\sigma)$ belongs to unstable (if $\tau=s$) or stable (if $\tau=u$) manifold of some node fixed point (sink or source, respectively), say $N$, and the topological closure of $Sep^{\tau}(\sigma)$ is a topologically embedded $d$-sphere that equals $W^{\tau}(\sigma)\cup\{N\}$.
\end{prop}

For $1\le m\le n$, we presume Euclidean space $\mathbb{R}^m$ to be included naturally in $\mathbb{R}^n$ as the subset whose final $(n-m)$ coordinates each equals $0$. Let $e: M^m\to N^n$ be an embedding of closed $m$-manifold $M^m$ in the interior of
$n$-manifold $N^n$. One says that $e(M^m)$ is \textit{locally flat at} $e(x)$, $x\in M^m$, if there exists a neighborhood $U(e(x))=U$ and a homeomorphism $h: U\to\mathbb{R}^n$ such that $ h(U\cap e(M^m))=\mathbb{R}^m\subset\mathbb{R}^n. $
Otherwise, $e(M^m)$ is \textit{wild at} $e(x)$ \cite{DavermanVenema-book-2009}. The similar notation for a compact $M^m$, in particular $M^m=[0;1]$).

Note that a separatrix $Sep^{\tau}(\sigma)$ is a smooth manifold. Hence, $Sep^{\tau}(\sigma)$ is locally flat at every point \cite{DavermanVenema-book-2009}. However a-priori, $clos~Sep^{\tau}(\sigma)=W^{\tau}(\sigma)\cup\{N\}$ could be wild at the unique point $N$.

One of the key statement for proving Theorem \ref{thm:rough-decomposition} is the following result proved in \cite{BonattiGrinesMedvedevPecou2002b} for $n=3$.
\begin{prop}\label{prop:nbhd-of-separatrix}
Let $e: \mathbb{S}^{n-1}\to M^n$ be a topological embedding of the $(n-1)$-sphere, $n\geq 3$, which is a smooth immersion everywhere, except at one point, and let $\Sigma^{n-1}=e(\mathbb{S}^{n-1})$. Then any neighborhood of $\Sigma^{n-1}$ contains a closed neighborhood of $\Sigma^{n-1}$ diffeomorphic to $\mathbb{S}^{n-1}\times [0;1]$.
\end{prop}
\textsl{Proof}. It is enough to prove the statement for $n\geq 4$. Let $\Sigma^{n-1}$ be a topologically embedded $(n-1)$-sphere that is smooth everywhere, except at one point, say $N\in\Sigma^{n-1}$. According \cite{Cantrell63} (see also \cite{CantrellEdwards65,Chernavskii66}), a wildly embedded $(n-1)$-sphere have to contain infinitely many points where the locally flatness fails provided $n\geq 4$. Therefore, $\Sigma^{n-1}$ is a locally flat embedded $(n-1)$-sphere. This completes the proof.
$\Box$

The following statement proved in \cite{GrinesMedvedevPochinkaZh2010} gives the sufficient condition for a Morse-Smale diffeo\-mor\-p\-hism to be polar.
\begin{prop}\label{prop:suff-condition-for-polar}
Let $f: M^n\to M^n$ be a Morse-Smale diffeomorphism without codimension one saddle periodic orbits. Then $f$ is a polar diffeomorphism, i.e. $f$ has a unique source periodic orbit and unique sink periodic orbit. Moreover, $M^n$ is orientable.
\end{prop}

The following propositions was proved in \cite{MedUman1979} (Lemma 6 and Lemma 7). For the Reader convenience, we give a sketch of the proof.
\begin{prop}\label{prop:skew-from-annulus}
Let $B^n_1$, $B^n_2\subset S^n$ be disjoint $n$-ball such that their boundaries $S^{n-1}_1=\partial B^n_1$, $S^{n-1}_2=\partial B^n_2$ are locally flat embedded $(n-1)$-spheres. Then given any homeomorphism $\psi: S^{n-1}_1\to S^{n-1}_2$, the manifold $N^n$ obtained from $S^n\setminus (B^n_1\cup B^n_2)$ after the identification of $S^{n-1}_1$, $S^{n-1}_2$ under $\psi$ is homeomorphic to
$S^{n-1}\otimes S^1$.
\end{prop} Since $S^{n-1}_1$, $S^{n-1}_2$ are locally flat embedded $(n-1)$-spheres, $S^n\setminus (B^n_1\cup B^n_2)$ is a closed $n$-annulus homeomorphic to $S^{n-1}\times S^1$ \cite{Wall-1969}. Hence, $N^n$ is homeomorphic to $S^{n-1}\otimes S^1$.
$\Box$

\begin{prop} \label{prop:connection-sum}
Let $M^n$ be a closed (topological) manifold and $S^{n-1}$ an $(n-1)$-sphere topologically imbedded in $M^n$. Suppose that $S^{n-1}$ has an open neighborhood $U$ homeomorphic to the prime product $S^{n-1}\times (-1,1)$. If the manifold $M^n\setminus U$ is connected, then there is the topological closed manifold $M^n_1$ such that $M^n$ is homeomorphic to the connected sum
$M^n= M^n_1\sharp (S^{n-1}\otimes S^1)$.
\end{prop}
\textsl{Proof}. Denote by $S^{n-1}_1$, $S^{n-1}_2$ the components of the boundary $\partial U$. Clearly that $S^{n-1}_1$, $S^{n-1}_2$ are $(n-1)$-spheres locally flat embedded in $M^n$. Since $M^n\setminus U$ is connected, there is a closed subset $D^n$ containing $S^{n-1}_1$, $S^{n-1}_2$ such that $D^n$ is homeomorphic to $clos~B^n\setminus (B^n_1\cup B^n_2)$ where $B^n_1$, $B^n_2\subset B^n$ are disjoint $n$-balls. Here, $S^{n-1}_i=\partial B^n_i$, $i=1,2$. We see that the boundary $\partial D^n$ consists of three components $S^{n-1}_0$, $S^{n-1}_1$, $S^{n-1}_2$ each homeomorphic an $(n-1)$-sphere. Attaching a closed $n$-ball to $D^n$ along the component $S^{n-1}_0$, one gets the set homeomorphic $S^n\setminus (B^n_1\cup B^n_2)$. The result now follows from Proposition \ref{prop:skew-from-annulus}.
$\Box$


To prove Theorems \ref{thm:new-proj-decomposition-for-diff}, \ref{thm:new-proj-decomposition-for-flow} we need the following description of projective-like manifold.
\begin{prop}\label{prop:def-proj-like}
A closed manifold $M^n$, $n\geq 3$, is projective-like if and only if $M^n$ is a disjoin union of an open $n$-ball $B^n$ and $k$-sphere $S^k$, $1\leq k\leq n-1$, locally flat embedded in $M^n$.
\end{prop}
\textsl{Proof}. Let $M^n$ be a projective-like manifold. Thus, $n\in\{4,8,16\}$, and $M^n$ is a manifold obtained from the closed $n$-ball $clos~\mathbb{B}^{n}=\mathbb{B}^{n}\cup\partial\mathbb{B}^{n}$ after the identification of every fiber of the
Hoph fiber bundle $(\partial\mathbb{B}^{n}=S^{n-1},S^{\frac{n}{2}}, S^{\frac{n}{2}-1})$ with a point. Denote by $S$ the set obtained after this identification. Since the Hoph fiber bundle is locally trivial, $S$ is homeomorphic to the base $S^{\frac{n}{2}}$ that is an $\frac{n}{2}$-sphere. Clearly that $\partial\mathbb{B}^{n}=S^{n-1}$ is the locally flat embedded $(n-1)$-sphere in the closed $n$-ball $clos~\mathbb{B}^{n}$. Hence, $S$ is also locally flat embedded in $M^n$. We see that $M^n$ is the disjoint union of the open $n$-ball $B^n$ and $\frac{n}{2}$-sphere $S$ that is locally flat embedded in $M^n$.

Now, suppose that $M^n$ is a disjoin union of an open $n$-ball $B^n$ and $k$-sphere $\Sigma^k$, $1\leq k\leq n-1$, locally flat embedded in $M^n$. Then $\Sigma^k$ has an open tubular neighborhood $T(\Sigma^k)$ such that its boundary $\partial T(\Sigma^k)$ is a submanifold of codimension one, and $T(\Sigma^k)$ is the total space of a locally trivial fiber bundle with the base $\Sigma^k$ and a fiber $B^{n-k}$ \cite{Hirsch-book-1976}. For convenience, we can assume that each fiber $B^{n-k}$ is an $(n-k)$-ball such that the boundary $\partial B^{n-k}=S^{n-k-1}$ belongs to $\partial T(\Sigma^k)$, and the center of $B^{n-k}$ belongs to $\Sigma^k$.

First, we have to show that $\partial T(\Sigma^k)$ is homeomorphic to $S^{n-1}$.
Let us construct flows $f^t_0$ and $f^t_1$ on the sets $B^n$ and $clos~T(\Sigma^k)=T(\Sigma^k)\cup\partial T(\Sigma^k)$ respectively, as follows. Take an arbitrary point $x_0\in B^n$ that does not belong to $clos~T(\Sigma^k)$. Since $B^n$ is an open ball, there is a flow $f^t_0$ on $B^n$ such that $f^t-0$ has a unique fixed point $x_0$ that is a source, and all one-dimensional trajectories leave any compact part of $B^n$ in the positive direction (time increases).
The flow $f^t_1$ is arranged as follows: a) each disk $\tilde{B}^{n-k}$ that is a fiber of the locally trivial bundle $\left(T(\Sigma^k),\Sigma^k,B^{n-k}\right)$) is invariant under $f^t_1$; b) the restriction of $f^t_1$ on $\tilde{B}^{n-k}$ has a sink at a point on $\Sigma^k$, corresponding to the center of the disk $B^{n-k}$, and has the set of equilibria that fill out the entire boundary of the disc $\tilde{B}^{n-k}$; c) the one-dimensional trajectories on the set $(T(\Sigma^k)\setminus\Sigma^k)\cap\tilde{B}^{n-k}$ move in the positive direction to the sinks.

Let $\tilde{\Sigma}^{n-1}$ be an $(n-1)$-sphere locally flat imbedded in $M^n$ such that $\tilde{\Sigma}^{n-1}$ bounds $n$-ball $b_0^n$ with a point $x_0$ inside, and $\tilde{\Sigma}^{n-1}$ is transversal (in the topological sense) to the trajectories of the flow $f^t_0$. From the properties of this flow, and the equality $M^n=\Sigma^k\cup B^n$, and from the fact that $\partial T(S^k)$ is a compact subset of $B^n$, it follows that there is the number $\tau>0$ such that the set $\tilde{\Sigma}^{n-1}_\tau=f^\tau_0(\tilde{\Sigma}^{n-1})$
belongs to $T(\Sigma^k)$. Moreover, the set $\partial T(\Sigma^k)$ belongs to $f^\tau_0(b_0^n)$. Clearly, $\tilde{\Sigma}^{n-1}_\tau$
is an $(n-1)$-sphere that locally flat embedded in $M^n$. We can assume that $\tilde{\Sigma}^{n-1}_\tau$ belongs to the wandering set of the flow $f^t_1$.

The intersection  $T(\Sigma^k)\cap f^\tau_0(b_0^n)$ is an open set whose boundary contains $\partial T(\Sigma^k)$. Since
$\partial T(\Sigma^k)$ is a submanifold of codimension one, $\partial T(\Sigma^k)$ has a semi-neighborhood
$\mathfrak{U}\subset\left(T(\Sigma^k)\cup\partial T(\Sigma^k)\right)\cap f^\tau_0(b_0^n)$ in the set
$clos\,\left(T(\Sigma^k)\cap f^\tau_0(b_0^n)\right)$ that is homeomorphic to
$(0;1]\times\partial T(\Sigma^k)$. Let us take an open subset $int\,\mathfrak{U}\subset\mathfrak{U}$ homeomorphic to
$(0;1)\times\partial T(\Sigma^k)$. Obviously,

\begin{equation}\label{eq:fundamental-groups-neghb-boundary}
   \pi_i(int\,\mathfrak{U})=\pi_i\left((0;1)\times\partial T(\Sigma^k)\right)=\pi_i\left(\partial T(\Sigma^k)\right),
   \quad i=0,\ldots, n-2.
\end{equation}

The set $A_0=B^n\setminus clos\,f^\tau_0(b_0^n)$ is an open $n$-dimensional annulus homeomorphic to $(0;1)\times\mathbb{S}^{n-1}$. Therefore, its homotopy groups $\pi_i(A_0)$ are equal to zero for all $i=0,\ldots,n-2$. Let us consider a representative
$\gamma: S^i\to int\,\mathfrak{U}$ of the group $\pi_i(int\,\mathfrak{U})$ where $S^i$ is an $i$-sphere. Since
$\gamma(S^i)\cap\partial T(\Sigma^k)=\emptyset$, there exists a number $\tau_1>0$ such that $f^{\tau_1}(\gamma(S^i))\subset A_0$. Therefore, (\ref{eq:fundamental-groups-neghb-boundary}) implies that $\pi_i(\partial T(\Sigma^k))=0$ for all $i=0,\ldots,n-2$. It follows from the validity of Poincare conjecture for all dimensions $n\geq 3$ (see \cite{Freedman1982,Newman1966,Perelman2003a,Perelman2003b,Smale61}) that the set $\partial T(\Sigma^k)$ is homeomorphic to an $(n-1)$-sphere.

Since $\partial T(\Sigma^k)$ is homeomorphic to $S^{n-1}$, the locally trivial bundle $(T(\Sigma^k),\Sigma^k,B^{n-k})$ is (globally) non-trivial. The projection $\pi : T(\Sigma^k)\to \Sigma^k$ of this bundle induces the projection
$\pi_*: \partial T(\Sigma^k)\to\Sigma^k$ such that $\pi_*^{-1}(x)=\partial B^{n-k}=S^{n-k-1}$ for any $x\in\Sigma^k$. Since the bundle
$(T(\Sigma^k),\Sigma^k,B^{n-k})$ is locally trivial, $\pi_*$ induces the locally trivial bundle
$(\partial T(\Sigma^k),\Sigma^k,\partial B^{n-k})=(S^{n-1},S^k,S^{n-k-1})$.
According \cite{Adams1960} (see also \cite{Novikov-book-1976}), there are only following such bundles
 $$ S^3\to S^2, \mbox{ fiber } S^1;\qquad S^7\to S^4, \mbox{ fiber } S^3;\qquad S^{15}\to S^8, \mbox{ fiber } S^7. $$
It is easy to see that these bundles correspond to the following pairs $(n,k)$: $(4,2)$, $(8,4)$, $(16,8)$.

We see that $M^n$ is a disjoin union of an open $n$-ball $B^n$ and $\frac{n}{2}$-sphere $S^{\frac{n}{2}}$ locally flat embedded in $M^n$ where $n\in\{4,8,16\}$. One can consider the spheres $S^{\frac{n}{2}}$, $n=4,8,16$ being bases of the Hoph bundles $(S^{n-1},S^{\frac{n}{2}}, S^{\frac{n}{2}-1})$. Since $\partial\mathbb{B}^{n}=S^{n-1}$, $M^n$ can obtained from the closed $n$-ball $clos~\mathbb{B}^{n}=\mathbb{B}^{n}\cup\partial\mathbb{B}^{n}$ after the identification of every fiber of the
Hoph bundle $(\partial\mathbb{B}^{n}=S^{n-1},S^{\frac{n}{2}}, S^{\frac{n}{2}-1})$ with a point. This completes the proof.
$\Box$

\section{Proofs of the main results}\label{s:proofs}

\begin{lm}\label{lm:exist-saddle-with-no-hetero}
Let $f: M^n\to M^n$ be a Morse-Smale diffeomorphism without heteroclinic submanifolds on codimension one separatrices. Suppose $\sigma$ is a codimension one saddle periodic point such that $\dim W^s(\sigma)=1$, $\dim W^u(\sigma)=n-1$. Then there exists a codimension one saddle periodic point $\sigma_*$ such that $Sep^u(\sigma_*)$ has no heteroclinic intersections.
\end{lm}
\textsl{Proof}. Without loss of generality, one can assume that all periodic points are fixed. Given any $p$, $q\in NW(f)$, we put $p\prec q$ provided $W^s(p)\cap W^u(q)\neq\emptyset$ and there no other $r\in NW(f)$ such that $W^s(p)\cap W^u(r)\neq\emptyset$, $W^s(r)\cap W^u(q)\neq\emptyset$. One knows that $\prec$ is a partial ordering and this ordering is strict \cite{Smale60a,Smale67}.

Suppose that $Sep^u(\sigma)$ has heteroclinic intersections (otherwise, nothing to prove). The chain $\sigma\prec\sigma_1\prec\cdots$ has a maximum point, say $\sigma_*$. Since $f$ has no heteroclinic manifolds on codimension one separatrices, every saddle in this chain has a codimension one unstable separatrix. Since $\sigma_*$ is a maximum point in the chain above, $Sep^u(\sigma_*)$ has no heteroclinic intersections.
$\Box$

\medskip
\textsl{Proof of Theorem \ref{thm:rough-decomposition}}. Taking a sufficiently large iteration, if necessary, one can assume that every periodic point of $f$ is fixed. If $\nu=0$ then $f$ is a polar Morse-Smale diffeomorphism by Proposition \ref{prop:suff-condition-for-polar}, and nothing to prove. Suppose now that $\nu\geq 1$.
Due to Lemma \ref{lm:exist-saddle-with-no-hetero}, there exists a codimension one saddle $\sigma$ whose codimension one separatrix has no heteroclinic intersections. For definiteness, let us assume that $\dim W^s(\sigma)=1$, $\dim W^u(\sigma)=n-1$. According to Proposition \ref{prop:GrinZhMedv2003b-top-closure-of-separatrices}, the topological closure $clos~W^u(\sigma)$ of $W^u(\sigma)$ is a topologically embedded $(n-1)$-sphere consisting of $W^u(\sigma)$ and a sink $\omega$. By Proposition \ref{prop:nbhd-of-separatrix}, there is an open neighborhood $U$ of $clos~W^u(\sigma)$ such that the topological closure $clos~U$ is diffeomorphic to
$\mathbb{S}^{n-1}\times [0;1]$. Since $clos~U$ contains the sink $\omega$, $f^k(clos~U)\subset U$ for a sufficiently large $k$. Passing to the iteration $f^k$ if necessary, we can assume, without loss of generality, that $k=1$.

Let us remove the neighborhood $U$ from the manifold $M^n$.
The manifold $M^n\setminus U$ has two boundary components $\Sigma^{n-1}_1$, $\Sigma^{n-1}_2$ each homeomorphic to $\mathbb{S}^{n-1}$.
Gluing to each $\Sigma^{n-1}_i$ an $n$-ball $B^n_i$, $i=1,2$, we get a smooth closed manifold $M^n_1$. Since $f(clos~U)\subset U$, we can extend the diffeomorphism $f$ to the manifold $M^n_1$ such that inside each ball $B^n_1$, $B^n_2$ the diffeomorphism we obtained $f_1: M^n_1 \to M^n_1$ has exactly one hyperbolic sink, while all points except the sinks are wandering.
Comparing the non-wandering sets of $f_1$ and $f$, one can see that $f_1$ has one saddle of codimension one less and one node (the sink, in this case) more. We call the described procedure a \textit{cutting along an unstable separatrix}. The similar cutting operation is considered along a stable separatrix (with adding sources).

After $\nu$ cuttings along codimension one separatrices of all saddles with the Morse indexes $n-1$ and $1$, we obtain a Morse-Smale diffeomorphism $f_{\nu}: M^n_{\nu}\to M^n_{\nu}$ of closed manifold $M^n_{\nu}$ consisting of finitely many connected components. The non-wandering set of $f_{\nu}$ contains exactly $\mu+\nu$ nodes, and does not contain codimension one saddles. By Proposition \ref{prop:suff-condition-for-polar}, each connected component of the manifold $M^n_{\nu}$ admits a polar diffeomorphism with exactly one source and exactly one sink. Hence, the number of connected components of $M^n_{\nu}$ is equal to $k=\frac{1}{2}\left(\mu + \nu\right)$. Therefore, the number $\mu+\nu$ is even.

The cutting procedure above allows to rebuild the original manifold $M^n$ from the obtained connected components. If one gets a connected manifold after the cutting along a codimension one separatrix, the intermediate manifold is homeomorphic to the connected sum of some closed manifold and $S^{n-1}\otimes S^1$, according to Proposition \ref{prop:connection-sum}.
If after the cutting along a codimension one separatrix one gets a manifold consisting of two connected components, then the intermediate manifold is homeomorphic to the connected sum of two closed manifolds.

Denote by $N^n_1$, $\ldots$, $N^n_{k}$ the connected components of the manifold $M^n_{\nu}$. The number of cuttings that does not increase the number of the connected components is equal to
the integer number
 $$ g=\nu-\frac{1}{2}\left(\mu + \nu\right)+1=\frac{1}{2}\left(\nu - \mu +2\right)\geq 0. $$
Each such cutting corresponds to the summand $S^{n-1}\otimes S^1$ in the resulting connected sum. Therefore, the connected sum of $k$ manifolds $N^n_1$, $\ldots$, $N^n_{k}$ and $g$ copies of $S^{n-1}\otimes S^1$ gives $M^n$.

For $g=0$, there are two possibilities:
a) all manifolds $N^n_1$, $\ldots$, $N^n_{k}$ are homeomorphic to the sphere $\mathbb{S}^n$ and hence, $M^n$ is homeomorphic to $\mathbb{S}^n$;
b) among $N^n_1$, $\ldots$, $N^n_{k}$, there exist exactly $1\leq l\leq k$ manifolds that are not homeomorphic to $\mathbb{S}^n$.

If $g\neq 0$, one can prove in a similar way that either $M^n$ is homeomorphic to the connected sum
$ \left(S^{n-1}\otimes S^1\right)\sharp\cdots\sharp\left(S^{n-1}\otimes S^1\right)$
of $g$ copies $S^{n-1}\otimes S^1$, or $M^n$ can be represented as the connected sum
$\left(S^{n-1}\otimes S^1\right)\sharp\cdots\sharp\left(S^{n-1}\otimes S^1\right)\sharp N^n_1\sharp\cdots\sharp N^n_{l}$
for some $1\leq l\leq k$, where each manifold $N^n_i$ ($i=1,\dots,l$) admits a polar diffeomorphism without codimension one saddle periodic orbits.

We remind Morse inequalities \cite{Smale60a}.
Let $M_j$ be the number of periodic points $p\in Per~(f)$ those stable Morse index equals $j = \dim W^s(p)$, and $\beta_i=rank~H_i(M^n,\mathbb{Z})$ the Betti numbers. Then
\begin{equation}\label{eq:Morse-Smale-inequalities1}
    M_0\ge \beta _0,\quad M_1 - M_0\ge\beta _1 - \beta _0,\quad M_2 - M_1 + M_0\ge \beta _2 - \beta _1 + \beta _0,\cdots
\end{equation}
\begin{equation}\label{eq:Morse-Smale-inequalities2}
    \sum _{i=0}^n(-1)^iM_i = \sum _{i=0}^n(-1)^i\beta _i.
\end{equation}
Since $N^n_i$ admits a polar Morse-Smale diffeomorphism without codimension one saddle periodic orbits, $M_1=0$. It follows from (\ref{eq:Morse-Smale-inequalities1}) that $\beta_0\leq M_0\leq\beta_0-\beta_1$.
Hence, $\beta_1(N^n_i) =0$. By duality, $\beta_{n-1}(N^n_i)=0$. This completes the proof.
$\Box$

\medskip
\textsl{Proof of Theorem \ref{thm:new-proj-decomposition-for-diff}}. Since the non-wandering set $NW(f)$ contains a single saddle fixed point, say $\sigma$, that is not codimension one, $l=1$ in the decompositions (\ref{eq:decomposition-1}), (\ref{eq:decomposition-3}) and the summand $N^n_1$ admits a polar diffeomorphism $f_1: N^n_1\to N^n_1$ the non-wandering set of whose consists of a sink $\omega_1$, a source $\alpha_1$, and the saddle fixed point $\sigma$. Moreover, $N^n_1\neq S^n$ because of $N^n_1$ contains two spheres $W^u(\sigma)\cup\{\omega_1\}$ and $W^s(\sigma)\cup\{\alpha_1\}$ transversally intersected at a unique point $\sigma$. Hence, the decomposition (\ref{eq:decomposition-2}) does not hold. In addition, due to Proposition \ref{prop:suff-condition-for-polar}, $N^n_1$ is orientable.
It follows from \cite{BonattiGrinesMedvedevPecou2002b} that if a Morse-Smale diffeomorphism of a closed 3-manifold has no heteroclinic intersections, then the number of periodic points can not be three and five. As a consequence, $\dim M^n\dim N^n_1=n\neq 3$. Later on, $n\geq 4$.

It follows from Proposition \ref{prop:GrinZhMedv2003b-top-closure-of-separatrices} that $\Sigma^k=W^u(\sigma)\cup\{\omega_1\}$ is a topologically embedded $k$-sphere where $2\leq k=\dim W^u(\sigma)\leq n-2$. It is well known \cite{grinmedpoch-book-2016-eng,Smale67} that a manifold is the disjoint union of unstable manifolds of non-wandering orbits. Hence, the manifold
 $$ N^n_1\stackrel{\rm def}{=}N^n=W^u(\sigma)\cup W^u(\omega_1)\cup W^u(\alpha_1)=W^u(\sigma)\cup\{\omega_1\}\cup W^u(\alpha_1)=\Sigma^k\cup W^u(\alpha_1) $$
is the disjoint union of the open $n$-ball $B^n=W^u(\alpha_1)$ and $k$-sphere $\Sigma^k$ topologically embedded in $N^n$.

The remaining assertions follows from \cite{MedvedevZhuzhoma2013-top-appl} (see also \cite{MedvedevZhuzhoma2012-jdcs}). For the Reader convenience, we give sketches of rest proofs.
Since $2\leq k=\dim W^u(\sigma)\leq n-2$, $M_0=M_n=M_k=1$. For $f^{-1}$, one holds $M_0=M_n=M_{n-k}=1$ and for $j\neq 0$, $n$, $k$, $n-k$, one holds $M_j=0$.
Since the left parts of (\ref{eq:Morse-Smale-inequalities2}) for $f$ and $f^{-1}$ are equal, $(-1)^k=(-1)^{n-k}$. Hence, $n=2m$ is even, where $m\ge 2$.
Let us show that $k=m$. Suppose the contradiction. Assume for definiteness that $k>m$. It follows from (\ref{eq:Morse-Smale-inequalities1}) that $\beta_1=\ldots=\beta_{n-k-1}=0$ because of $M_1=\ldots=M_{n-k-1}=0$. The Poincare duality implies that $\beta_1=\ldots=\beta_{k-1}=0$. Hence, $\beta_i=0$ for all $i=1$, $\ldots$, $n-1$. Then (\ref{eq:Morse-Smale-inequalities2}) becomes $1+(-1)^k+(-1)^n=1+(-1)^n$. This is impossible. Hence, $k=\frac{n}{2}$.

It remains to prove that the spheres $\Sigma^k=W^u(\sigma)\cup\{\omega_1\}=S^k$, $S^{n-k}=W^s(\sigma)\cup\{\alpha_1\}$ are locally flat provided $n\geq 6$.
It follows from \cite{Chernavskii66} (see \cite{CantrellEdwards65,Stallings63}) that $k$-manifold has no isolated wild points provided $n\ge 5$, $k\neq n-2$. As a consequence, $S^k$, $S^{n-k}$ are locally flat embedded $k$-spheres. The Theorem is proved.
$\Box$

\medskip
\textsl{Proof of Theorem \ref{thm:new-proj-decomposition-for-flow}}. Since the Morse-Smale flow $f^t$ has no periodic trajectories, the time one shift along the trajectories is a Morse-Smale diffeomorphism, say $f$. We keep the notation of the proof of Theorem \ref{thm:new-proj-decomposition-for-diff}. It was proved in \cite{MedvedevZhuzhoma2013-top-appl} (see also \cite{MedvedevZhuzhoma2012-jdcs}) that if $f$ is generated by a Morse-Smale flow, then the spheres $\Sigma^k=W^u(\sigma)\cup\{\omega_1\}=S^k$, $S^{n-k}=W^s(\sigma)\cup\{\alpha_1\}$ are locally flat provided $n\geq 4$. Now, the result follows from Theorem \ref{thm:new-proj-decomposition-for-diff}.
$\Box$

\medskip
\textsl{Proof of Theorem \ref{thm:realization}}. Previously, we introduce some notation. Let $M^n_1$, $M^n_2$ be $n$-manifolds supporting vector fields $v_1$, $v_2$ respectively. Suppose that $v_1$, $v_2$ are consistent on the boundaries $\partial M^n_1$, $\partial M^n_2$. For example, $v_1$ is outside and perpendicular to $\partial M^n_1$ while $v_2$ is inside and perpendicular to $\partial M^n_2$. Assume that there is a diffeomorphism
$\varphi: \partial M^n_1\to\partial M^n_2$. On the manifold $M^n_1\cup_{\varphi(\partial M^n_1)=\partial M^n_2}M^n_2$, the vector fields $v_1$, $v_2$ form the vector field denoted by $v_1\sharp v_2$. Below, all fixed points are hyperbolic. Saying that a vector field induces a diffeomorphism, we mean that the diffeomorphism is the shift-one-time along the trajectories of the flow induced by the vector field.

Let $V_{sink}$ be the vector field on the closed $n$-disk $D^n$ with a sink at the center such that $V_{sink}$ is inside and perpendicular to $\partial D^n=S^{n-1}$. Denote $-V_{sink}$ by $V_{source}$. Obviously, $V_{source}$ has a source at the center and $V_{source}$ is outside and perpendicular to $\partial D^n=S^{n-1}$. Then the vector field $V_{sink}\sharp V_{source}=V_{NS}$ on the $n$-sphere $S^n=D^n\cup_{id(S^{n-1})=S^{n-1}}D^n$ induces the Morse-Smale diffeomorphism $f_{SN}: S^n\to S^n$ of the North-Sough type.

First, one considers $g=0$. The diffeomorphism $f_{SN}$ is desired for the particular case $\nu=0$, $\mu=2$. For $\nu\geq 1$, one considers the vector field $V_{\nu,sink}$ on $D^n$ with $\nu$ codimension one saddles, one sink,
and $\nu$ sources such that $V_{\nu,sink}$ is inside and perpendicular to $\partial D^n=S^{n-1}$. Each saddle has two stable one-dimensional separatrices and an $(n-1)$-dimensional unstable separatrix that together with the sink surrounds one source. In addition, one can assume that codimension one separatrices have no heteroclinic intersections.
The vector field $V_{source}\sharp V_{\nu,sink}$ induces the desired Morse-Smale diffeomorphism $S^n\to S^n$ with $\nu$ codimension one saddles and $\nu+2=\mu$ nodes.

Denote by $V_{a}$ the vector field on $D^{n-1}\times S^1$ that is inside and perpendicular to the boundary
$\partial (D^{n-1}\times S^1)=S^{n-2}\times S^1$ such that $V_{a}$ has $a\geq 1$ sinks and $a$ codimension one saddles. In addition, each saddle has an $(n-1)$-dimensional stable separatrix and two one-dimensional unstable separatrices going to sinks, and the stable separatrix intersects $S^{n-2}\times S^1$ at $S^{n-2}\times\{x\}$ for some $x\in S^1$. Take a copy $D^{n-1}\times S^1$ supporting the vector field $-V_{b}$. Clearly, $-V_{b}$
is outside and perpendicular to the boundary $\partial (D^{n-1}\times S^1)=S^{n-2}\times S^1$ and $V_{b}$ has $b\geq 1$ sources and $b$ codimension one saddles. Denote by $R_{\beta}: S^1\to S^1$ the rigid rotation
$x\to x+\beta\, mod\, 1$ which induces the diffeomorphism $\psi_{\beta}: S^{n-2}\times S^1\to S^{n-2}\times S^1$ as follows $\psi_{\beta}(z,x)=(z,R_{\beta}(x))$. On the closed manifold
$$ (D^{n-1}\times S^1)\bigcup_{\psi_{\beta}(\partial (D^{n-1}\times S^1))}(D^{n-1}\times S^1)=S^{n-1}\times S^1, $$
the vector fields $V_{a}$, $-V_{b}$ form the vector field $W_{a+b}$ with $a+b$ codimension one saddles and $a+b$ nodes. One can choose $\beta$ such that $W_{a+b}$ becomes a Morse-Smale vector field without intersections on codimension one separatrices. This $W_{a+b}$ induces the Morse-Smale diffeomorphism $f_{a,b}$ without intersections on codimension one separatrices. Actually, $f_{a,b}$ is gradient-like.

Now, one considers $g\geq 1$. In this case, $\nu=2g+(\mu-2)\geq 2$ because of the inequality $\mu\geq 2$. Note that since $g\in\mathbb{Z}$, the number $\nu-\mu\geq 0$ is even.
Denote by $[x]$ the integral part of $x\in\mathbb{R}$. By construction, the Morse-Smale diffeomorphism
$$ f_{\left[\frac{\nu}{2}\right],\left[\frac{\nu+1}{2}\right]}: S^{n-1}\times S^1\to S^{n-1}\times S^1 $$
has $\left[\frac{\nu}{2}\right] + \left[\frac{\nu+1}{2}\right]=\nu$ codimension one saddles and $\nu$ nodes consisting of $\left[\frac{\nu}{2}\right]$ sinks and $\left[\frac{\nu+1}{2}\right]$ sources. If $\nu=\mu$, $f_{\left[\frac{\nu}{2}\right],\left[\frac{\nu+1}{2}\right]}$ is the desired Morse-Smale diffeomorphism $S^n\to S^n$ with $\nu$ codimension one saddles and $\nu=\mu$ nodes. Now, suppose $\nu-\mu\geq 2$. Let us delete from $S^{n-1}\times S^1$ the sufficiently small neighborhoods of $\frac{1}{2}(\nu - \mu)$ sources and $\frac{1}{2}(\nu - \mu)$ sinks where each neighborhood is homeomorphic to an open $n$-ball $B^n$. Let us identify the boundary of every deleted neighborhood of source with the boundary of deleted neighborhood of a sink such that the diffeomorphism $f_{\left[\frac{\nu}{2}\right],\left[\frac{\nu+1}{2}\right]}$ induces the Morse-Smale diffeomorphism $f$ on the obtained closed manifold $M^n$. It follows from Proposition \ref{prop:skew-from-annulus} that $M^n$ is the connected sum of $1+\frac{\nu-\mu}{2}=g$ copies of $S^{n-1}\times S^1$. Calculations show that the non-wandering set $NW(f)$ consists of $\nu$ codimension one saddles and $\nu-2\cdot\frac{\nu-\mu}{2}=\mu$ nodes. This completes the proof.
$\Box$

\medskip
\textsl{Proof of corollaries \ref{cor:cor1-conditions-for-existing-heter-curves}, \ref{cor:periodic-trajectories -exist-if}}.
Outline of the proof of the corollaries is the same: if we assume the contrary, then there exist the decompositions of the ambient manifold $M^n$, according to Theorem \ref{thm:rough-decomposition}. It follows from Van Kampen Theorem (see, exm., \cite{Novikov-book-1976}), that $\pi_1(M^n)$ contains the subgroup $\mathbb{Z}\ast\cdots\ast\mathbb{Z}$. This contradiction proves the required assertions.
$\Box$

\bigskip
{\it E-mail:} vgrines@yandex.ru

\noindent
\textit{E-mail:} medvedev@unn.ac.ru

\noindent
\textit{E-mail:} zhuzhoma@mail.ru

\end{document}